\documentclass[12pt]{article}

\usepackage[T2A]{fontenc}
\usepackage[utf8]{inputenc}
\usepackage{graphics}
\usepackage{floatflt} 
\usepackage{color}
\usepackage{subfigure}
\usepackage{amssymb, amsmath, amsthm, bbm, verbatim} 
\newtheorem{teo}{Theorem}
\newtheorem{prop}{Proposition}
\newtheorem{df}{Definition}
\newtheorem{cor}{Corollary}
\newtheorem{rem}{Remark}
\newtheorem{lem}{Lemma}

\renewcommand{\r}{{\mathbb R}}
\newcommand{\z} {{\mathbb Z}}
\newcommand{\cn} {{\mathbb C}}

\newcommand{\be}{\begin{equation}}
\newcommand{\ee}{\end{equation}}
\newcommand{\bega}{\begin{gather}}
\newcommand{\enga}{\end{gather}}
\newcommand{\ban}{\begin{eqnarray*}}
\newcommand{\ean}{\end{eqnarray*}}

\newcommand{\supp}{\mbox{supp}\,}

\title{Uncertainty product for  Vilenkin groups
}
\author{
 Ivan Kovalyov 
\footnote{Department of Mathematics, National Pedagogical Dragomanov University,
Kiev, Pirogova 9, 01601, Ukraine},
Elena Lebedeva
\footnote{Faculty of Applied Mathematics and Control Processes, Saint Petersburg State University,
Universitetskaya nab., 7-9, Saint Petersburg,   199034, Russia; St. Petersburg Polytechnical University, Department of Calculus, 
 Polytekhnicheskay 29, 195251, St. Petersburg, Russia} 
}
\date{
i.m.kovalyov@gmail.com, ealebedeva2004@gmail.com}

\begin{document}
\maketitle

\begin{abstract}
We study a  localization of functions defined on  Vilenkin groups. To measure the localization we introduce two uncertainty products $UP_{\lambda}$ and $UP_{G}$ that are similar to the Heisenberg 
uncertainty product. $UP_{\lambda}$ and $UP_{G}$ differ from each other by the metric used for the Vilenkin group $G$. 
 We discuss  analogs of  a quantitative uncertainty principle. Representations for  $UP_{\lambda}$ and $UP_{G}$ in terms of Walsh and Haar basis are given.
\end{abstract}




\section{Introduction}
An uncertainty product for a function characterizes how concentrated is the function in time and frequency domain.    
Initially the notion of uncertainty product was introduced for $f\in L_2(\mathbb{R})$ by W.~Heisenberg \cite{H} and E.~Schr\"odinger \cite{Sch}. Later on extensions of this notion appeared for various algebraic
and topological structures. For periodic functions, it was suggested  by E.~Breitenberger  \cite{B}. For some particular cases of locally compact groups (namely a euclidean motion groups, non-compact semisimple Lie groups, Heisenberg groups) the counterpart was derived in \cite{PrSi}. 
Uncertainty products on compact Riemannian manifolds was discussed in \cite{Erb}.  In   
\cite{27},  this concept was introduced for functions defined on the Cantor group. In this paper, we discuss localization of functions defined on Vilenkin groups. 

To measure the localization we introduce a functional that is similar to the Heisenberg 
uncertainty product (see Definition \ref{5def1}). It depends  on the metric used for the Vilenkin group $G$. Two equivalent metrics are in common use for the group $G$. So we discuss two uncertainty products $UP_{\lambda}$ and $UP_{G}.$  The first one is a strict counterpart of ``dyadic uncertainty constant'' introduced in \cite{27} (see Theorems \ref{5UP} and \ref{5th1}). Usage of another metric in the second uncertainty product allows for exploitation of a modified Gibbs derivative that plays a role of usual derivative for the Heisenberg uncertainty product. At the same time it turns out that usage of Haar basis is a good approach for evaluation of  $UP_{G}$ (see Theorem \ref{UP_Haar}).  In particular, it allows for an estimate of Fourier-Haar coefficients for functions defined on the Vilenkin group (see Corollary \ref{Haar_coef}). The connection between 
$UP_{\lambda}$ and $UP_{G}$ is showed in Lemma \ref{comparing}.

\section{Auxiliary results}

We recall necessary facts about the Vilenkin group. More details can be found in \cite{GES,SWS}. 
The  Vilenkin group $G=G_p$, $p\in\mathbb{N}$, $p\neq 1$, is a set of  the sequences
%
$$
 x=(x_j)=(\dots,0,0,x_{-k},x_{-k+1},x_{-k+2},\dots),
$$
where $x_j\in \{0,\dots,p-1\}$ for $j\in \mathbb{Z}$.  The operation on $G$  is denoted by~$\oplus$ and defined as the
coordinatewise addition modulo $p:$
$$
(z_j)=(x_j)\oplus (y_j) \, \Longleftrightarrow\, z_j = x_j + y_j \,({\rm
mod}\, p)\quad\mbox{for}\quad j\in {\Bbb Z}.
$$
  The inverse operation of $\oplus$  is denoted by $\ominus$. 
The symbol $\ominus x$ denotes the inverse element of $x\in G.$  The sequence ${\bf 0}=(\dots,0,0,\dots)$ is a neutral element of $G$. If  $x\ne\bf 0$,  then there exists a unique number  $N=N(x)$ such that  $x_N\ne0$ and  $x_j=0$ for  $j<N$. The Vilenkin group 
 $G_p$, where $p=2$ is called the Cantor group. In this case the inverse operation  $\ominus$ coicides with the group operation $\oplus.$

Define a map  \, $\lambda:\,G\to [0,+\infty)$ 
$$
    \lambda(x)=\sum_{j\in \mathbb{Z}}x_jp^{-j-1}, \qquad x=(x_j)\in G.
$$
The mapping $x\mapsto \lambda(x)$ is a bijection taking $G\setminus\mathbb{Q}_0$ onto $[0,\,\infty),$ where
$\mathbb{Q}_0$ is a set  of all elements terminating with $p-1$'s.

 Two equivalent metrics are in common use for the group $G$. One metric is defined by $d_1(x,y):=\lambda(x\ominus y)$ for  $x,y\in G$. To define  another one  $d_2$ we consider a map  
 $\|\cdot\|_G:\ G\to [0,\infty)$, where  $\|{\bf 0}\|_G:=0$ and $\|x\|_G:=p^{-N(x)}$
  for $x\ne{\bf 0}$. Then $d_2(x,y):=\|x\ominus y\|_G$,  $x,y\in G$.    Given $n\in\z$ and $x\in G$,
denote by $I_{n}(x)$
the ball of radius $2^{-n}$ with the center at  $x$, i.e.
$$
I_{n}(x)=\{y\in G: d(x,y) < 2^{-n}\}.
$$
For brevity we set $I_j:=I_j({\bf 0})$ and $I:=I_0$.

We denote dilation on $G$ by $D:\ G\to G$,  and set $(Dx)_k=x_{k+1}$
for $x\in G$. Then $D^{-1}:\ G\to G$
 is the inverse mapping  $(D^{-1}x)_k=x_{k-1}$.  Set $D^k =
D\circ\dots\circ D$
($k$ times) if $k > 0$, and 
$D^k =
D^{-1}\circ\dots\circ D^{-1}$ ($-k$ times) if $k < 0$;
$D^0$ is the identity mapping.

We deal with functions taking $G$ to $\cn$. Denote
$\mathbbm{1}_E$  the characteristic function of a set  $E\subset G$.
Given  a function
$f:\ G\to \cn$ and a number $h\ge0$,   for every $x\in G$ we define
$f_{0, h}(x)=f(x\oplus \lambda^{-1}(h))$.
 Finally, we set  for $j\in\z$  
$$
f_{j,h}(x)=p^{j/2}f_{0,h}(D^jx), \quad x\in G.
$$

  The functional spaces
$L_q(G)$ and $L_q(E)$, where $E$ is a measurable subset of $G$,
are derived using the Haar measure (see~\cite{HR}). 

 Given $\xi\in G$, a group character of $G$ is defined by 
$$
\chi_\xi(x)=   \chi(x,\xi): = \exp\left(\frac{2\pi i}{p}\sum\limits_{j\in {\Bbb Z}} x_j\,\xi_{-1-j}\right).
$$
The functions ${\rm w}_n(x):=\chi(\lambda^{-1}(n),x)$ are called the  generalized Walsh functions. If $p=2$, than ${\rm w}_n$ are called the Walsh functions.
 
The Fourier transform of a function $f\in L^1(G)$ is defined by
\be \label{fourier}
   Ff(\omega)=\int\nolimits_Gf(x)\overline{\chi (x,\omega)}d\mu(x), \quad \omega\in G.
\ee
The Fourier transform is extended to ${ L}_2(G)$ in a
standard way, and the Plancherel equality takes place 
$$
\langle f,g\rangle:=
\int\limits_{G}f(x)\overline{g(x)}\,dx=
\int\limits_{G}F f(\xi)\overline{F g(\xi)}\,d\xi
=\langle F f, F g\rangle,\quad f,g\in L_2(G).
$$
The inversion formula is valid for any $f\in L_2(G)$
\begin{equation*}
F^{-1}F f(x) = \int_G  Ff(\omega) \chi (x,\,\omega) d \mu (\omega) = f(x).
\end{equation*}

It is straightforward to see that 
\begin{equation}
\label{6hat_jk}
	 F(f_{j,n})(\xi)= p^{-j/2} \chi(k,\,D^{-j}\xi)
	 Ff(D^{-j} \xi), \quad n\in \z_+, j\in\z.
\end{equation}

The discrete Vilenkin-Chrestenson transform of a vector
 $x=(x_k)_{k=\overline{0,p^n-1}}\in \cn^{p^n}$ is a vector $y=(y_k)_{k=\overline{0,p^n-1}}\in \cn^{p^n}$, where
\be
\label{71.1}
y_{k}
=p^{-n}\sum_{s=0}^{p^n-1} x_s {\rm w}_{k}(\lambda^{-1}(s/p^n)),
\quad 0\le k \le p^n-1.
\ee
The inverse transform is 
\be
\label{71.2}
x_{k}
=\sum_{s=0}^{p^n-1} y_s \overline{{\rm w}_{k}(\lambda^{-1}(s/p^n))}.
\quad 0\le k \le p^n-1.
\ee

Given $f:G_2\to\cn$, the function
$$
f^{[1]}(x):=\lim_{n\to \infty} \sum_{j=-n}^{n} 2^{j-1}(f(x)-f_{0, 2^{-j-1}}(x))
$$
is called the Gibbs derivative of a function $f$. The following properties hold true 
\be
\label{cl_gibbs}
Ff^{[1]}(\xi)=\lambda(\xi) Ff(\xi), \quad {\rm w}_n^{[1]}(x) = n {\rm w}_n(x). 
\ee

Set $\varphi=\mathbbm{1}_{I}$. The Haar functions $\psi^{\nu},$ $\nu=1,\dots,p-1$ are defined by 
\be
\label{haar}
\psi^{\nu}(x) = \sum_{n=0}^{p-1} {\rm exp}\left(\frac{2\pi i \nu n}{p}\right) \varphi(Dx\oplus \lambda^{-1}(n)).
\ee
The system $\psi_{j,k}^{\nu},$ $\nu=1,\dots,p-1,$ $j\in\z,$ $k\in\z_+,$ forms an orthonormal basis (Haar basis) for $L_2(G)$, see \cite{F,Lang}. 

It follows from (\ref{fourier}) that $F\varphi = \varphi = \mathbbm{1}_{I}$ and 
$F\psi=\mathbbm{1}_{I_{0}\oplus\lambda^{-1}(p-\nu)}.$ Taking into account (\ref{6hat_jk}), we get
\be
\label{hat_haar}
F\psi^{\nu}_{j,k}(\xi)= p^{-j/2} \chi(k,\,D^{-j}\xi) \mathbbm{1}_{I_{-j}\oplus\lambda^{-1}((p-\nu)p^j)}. 
\ee
Given $f\in L_1(G),$ the modified Gibbs derivative ${\cal D}$ is defined by 
\be
\label{gibbs}
F{\cal D}f = \|\cdot\|_G Ff.
\ee
It was introduced in \cite{G} for $L_1(G_2)$. Such kind of operators are often called pseudo-differential.

\begin{prop}
\label{eigen}
Suppose $g,$ $Fg,$ $\|\cdot\|_G Fg$ are locally integrable on $G$, $j\in\z$. Then the assertion
${\rm supp}\, \widehat{g} \subset I_{-j-1}\setminus I_{-j}$
is necessary and
sufficient for
$g$ to be an eigenfunction of ${\cal D}$
corresponding to the eigenvalue~$p^{j }$.
\end{prop}

The proof can be rewritten from Proposition 1 \cite{LS15}, where it is proved for the Cantor group.

\begin{cor}
\label{eigen_Haar}
Any Haar function $\psi_{j,k}^{\nu}$ is an eigenfunction of ${\cal D}^{\alpha}$ corresponding to the eigenvalue $p^{j}.$
\end{cor}

\textbf{Proof.} The statement follows from Proposition \ref{eigen} and (\ref{hat_haar}).\hfill $\Box$

\section{Uncertainty product and metrics}
Originally, the concept of an uncertainty product was introduced for the real line case in 1927. 
The Heisenberg uncertainty product  of $f \in L_2(\mathbb{R})$ is the functional 
$UC_H(f):=\Delta_{f}\Delta_{\widehat{f}}$ such that
$$ 
\Delta_{f}^2:=\|f\|^{-2}_{L^2(\mathbb{R})}\int_{\mathbb{R}}(x-x_{f})^2|f(x)|^2\,d x, \quad
\Delta_{\widehat{f}}^2:=\|\widehat{f}\|^{-2}_{L^2(\mathbb{R})}
 \int_{\mathbb{R}}(t-t_{\widehat{f}})^2|\widehat{f}(t)|^2\,d t, 
$$
$$
x_{f}:=\|f\|^{-2}_{L^2(\mathbb{R})}\int_{\mathbb{R}}x|f(x)|^2\,d x, \quad
t_{\widehat{f}}:=\|\widehat{f}\|^{-2}_{L^2(\mathbb{R})}\int_{\mathbb{R}}
t|\widehat{f}(t)|^2\,d t, 
$$  
where $\widehat{f}$ denotes the Fourier transform of $f\in L_2(\r).$
It is well known that 
$UC_H(f)\geq 1/2$ for a function $f \in L_2(\mathbb{R})$ and the minimum is attained on the Gaussian.  To motivate the definition of a localization characteristic for the Vilenkin group we note that on one hand  
 $x_f$ is the solution of the minimization problem
$$
\min_{\tilde{x}}\int_{\mathbb{R}}(x-\tilde{x})^2|f(x)|^2\,d x,
$$  
and on another hand 
the sense of the sign ``-'' in the definition of $\Delta_{f}$ is the distance between $x$ and $x_f.$
So we come to the main definition.

\begin{df}
\label{5def1}
Suppose $f : G \to \mathbb{C}$, $f\in L_2(G)$, and $d$ is a metric on $G$, then a functional   
$$
UP(f):=V(f)V(Ff), \quad \mbox{ where }
$$
$$
V(f):=
\frac{1}{\|f\|^2_{L_2(G)}} 
\min_{\tilde{x}}\int_{G}(d(x,\, \tilde{x}))^2|f(x)|^2\,d x
$$
is called the uncertainty product of a function $f$ defined on the Vilenkin group.
\end{df}

Thus, we study two uncertainty products $UP_{\lambda}$ and $UP_G$ that corresponds to the metric $d_1(x,y):=\lambda(x\ominus y)$ and $d_2(x,y):=\|x\ominus y\|_G$. More precisely, 
$$
UP_{\lambda}(f):=V_{\lambda}(f)V_{\lambda}(Ff), \quad \mbox{ where }
$$
$$
V_{\lambda}(f):=
\frac{1}{\|f\|^2_{L_2(G)}} 
\min_{\tilde{x}}\int_{G}(\lambda(x\ominus \tilde{x}))^2|f(x)|^2\,d x.
$$

The functional $UP_G$ is defined as
$$
UP_{G}(f):=V_{G}(f)V_{G}(Ff), \quad \mbox{ where }
$$
$$
V_{G}(f):=
\frac{1}{\|f\|^2_{L_2(G)}} 
\min_{\tilde{x}}\int_{G} \|x\ominus \tilde{x}\|_G^2|f(x)|^2\,d x.
$$

The functional $UP_{\lambda}$ for functions defined on the Cantor group was introduced and studied in \cite{27}. The following results are extended from the Cantor group to the Vilenkin group without any essential changes. So we omit the proofs. 
\begin{teo}
\label{5UP}
Suppose  $f: G \to \mathbb{C},$ 
$f\in L_2(G).$ Then the following inequality holds true 
$$
UP_{\lambda}(f)\geq C,
\mbox{ where }
C \simeq 8.5 \times 10^{-5}.
$$
\end{teo}  

\begin{teo}
\label{5th1}
Let
$
f(x)=\mathbbm{1}_{\lambda^{-1}[0,\,1)}(x) \sum_{k=0}^{\infty}a_k {\rm w}_k( x)
$
 be a uniformly convergent series. Denote
 $$
f_n(x)=\mathbbm{1}_{\lambda^{-1}[0,\,1)}(x) \sum_{k=0}^{p^n-1}a_k {\rm w}_k( x).
$$
 Let $V_{\lambda}(f)<+\infty, $ $V_{\lambda}(F{f})<+\infty.$ Then  
$
UP_{\lambda}(f)\,=\,\lim_{n \to \infty}V_{\lambda}(f_n)V_{\lambda}(F{f}_n),
$  where
$$
V_{\lambda}(f_n)=\frac{\min_{k_0=\overline{0,p^n-1}}\sum_{k=0}^{p^n-1} p^{-n} |b_{\lambda(\lambda^{-1}(k)\oplus \lambda^{-1}(k_0))}|^2 ((k+1)^3-k^3)/3}{\sum_{k=0}^{p^n-1}|a_k|^2},
$$
\begin{equation*}
\label{5finiteV}
V_{\lambda}(F{f}_n)=\frac{\min_{k_1=\overline{0,p^n-1}}\sum_{k=0}^{p^n-1} |a_{\lambda(\lambda^{-1}(k)\oplus \lambda^{-1}(k_1))}|^2 ((k+1)^3-k^3)/3}{\sum_{k=0}^{p^n-1}|a_{k}|^2},
\end{equation*}
and 
$b_{k}$,
 $0\le k \le p^n-1$,
is the inverse discrete Vilenkin-Chrestenson transform (\ref{71.2}).  
\end{teo}

The following Lemma shows that the functionals $UP_{\lambda}$ and $UP_{G}$ have the same order.

\begin{lem} \label{comparing}
Suppose $f\in L_2(G)$, then $p^{-4} UP_G(f) \le UP_{\lambda}(f)<  UP_G(f).$
\end{lem}
 
 \noindent
\textbf{Proof.} It is sufficient to note that $p^{-1} \|x\|_{G} \le \lambda (x) < \|x\|_{G}.$ \hfill $\Box$

Taking into account Theorem \ref{5UP}, we conclude that $UP_G$ has a positive lower   bound. So, $UP_G$ satisfies the uncertainty principle.


\noindent
{\bf Example 1}. Let us illustrate a definition of $UP_{G}$ for $p=2$ using functions $f_1,$ $g_1$, $f_2$, and $g_2$ taken from  \cite[Example 1]{27}. Recall $f_1(x)= \mathbbm{1}_{\lambda^{-1}[0,\,1/4)}(x),$ $g_1(x)= \mathbbm{1}_{\lambda^{-1}[3/4,\,1)}(x),$ $f_2(x)=\mathbbm{1}_{\lambda^{-1}[0,\,3/8)}(x),$ and $g_2(x)= \mathbbm{1}_{\lambda^{-1}[3/4,\,9/8)}(x).$ Their Walsh-Fourier transforms are
$Ff_1=\mathbbm{1}_{\lambda^{-1}[0,\,4)}/4,$
$F{g_1}={\rm w}_3\left(\cdot/4\right)\mathbbm{1}_{\lambda^{-1}[0,\,4)}/4,$
$F{f_2}=\mathbbm{1}_{\lambda^{-1}[0,\,4)}/4+{\rm w}_1\left(\cdot/4\right)\mathbbm{1}_{\lambda^{-1}[0,\,8)}/8,$ 
 and
 $F{g_2}={\rm w}_3\left(\cdot/4\right)\mathbbm{1}_{\lambda^{-1}[0,\,4)}/4+{\rm w}_1(\cdot)\mathbbm{1}_{\lambda^{-1}[0,\,8)}/8.$
Given $\alpha \in [0,\,\infty)$,
since the mapping $\alpha \mapsto \|\lambda^{-1}(\alpha)\|_{G}$ is increasing and a  measure of the set 
$\lambda^{-1}[a,\,b)\ominus \tilde{x}$ does not depend on $\tilde{x},$ 
it follows that  
$$
\min_{\tilde{x}}\int_{\lambda^{-1}[0,\,\frac14)} \|x\ominus \tilde{x}\|_G\,d x = 
\min_{\tilde{x}}\int_{\lambda^{-1}[0,\,\frac14)\ominus \tilde{x}} \|\tau\|_G\,d \tau =
\int_{\lambda^{-1}[0,\,\frac14)} \|\tau\|_G\,d \tau,
$$
and $\lambda^{-1}[0,\,1/4)$ is a set of minimizing $\tilde{x}$'s as well.
So, taking into account   
$\|f_1\|^2_{L_2(G)}=\|Ff_1\|^2_{L_2(G)}=1/4,$ we get
$$
 V_{G}(f_1)=\frac{1}{\|f_1\|^2_{L_2(G)}}\min_{\tilde{x}} \int_{G}\|x\ominus\tilde{x}\|^2_{G}|f_1(x)|^2dx=
4 \min_{\widetilde{x}} \int_{\lambda^{-1}[0,\frac{1}{4})}\|x\ominus\tilde{x}\|^2_{G}dx
$$
$$
=4\int_{\lambda^{-1}[0,\frac{1}{4})}\|\tau\|^2_{G}d\tau=
4\sum_{i=2}^{\infty}\int_{\lambda^{-1}[\frac{1}{2^{i+1}},\frac{1}{2^i})}\|\tau\|^2_{G}d\tau=4\sum_{i=2}^{\infty}\left(\frac{1}{2^i}-\frac{1}{2^{i+1}}\right)2^{-2i}=\frac{1}{28}.
$$ 
Analogously, we obtain
$$
 V_{G}(Ff_1)=\frac{1}{\|Ff_1\|^2_{L_2(G)}}\min_{\tilde{x}} \int_{G}\|x\ominus\tilde{x}\|^2_{G}|Ff_1(x)|^2dx=\frac{1}{4}
\min_{\tilde{x}} \int_{\lambda^{-1}[0,4)}\|x\ominus\widetilde{x}\|^2_{G}dx
$$
$$
=\frac{1}{4}\int_{\lambda^{-1}[0,4)}\|\tau\|^2_{G}d\tau =
\frac{1}{4}\sum_{i=-2}^{\infty}\int_{\lambda^{-1}[\frac{1}{2^{i+1}},\frac{1}{2^i})}\|\tau\|^2_{G}d\tau=\frac{1}{4}\sum_{i=-2}^{\infty}\left(\frac{1}{2^i}-\frac{1}{2^{i+1}}\right)2^{-2i}=\frac{64}{7}.
$$
Thus, $UP_G(f_1)=16/49.$
Using the same arguments, we calculate $UP_G$ for the remaining functions. We collect all the information in  Table \ref{table1}. Values of $UP_{\lambda}$ we extract from \cite[Example 1]{27}.
Columns named $\tilde{x}_0(f)$ and $\tilde{t}_0(f)$ contain  sets of $\tilde{x}$ and $\tilde{t}$ minimizing the functionals $V_{\lambda}(f)$, $V_G(f)$ and $V_{\lambda}(Ff)$, $V_G(Ff)$ respectively. With respect both uncertainty products $UP_G$ and $UP_{\lambda}$, functions $f_1$ and $g_1$ have the same localization, while function $f_2$ is more localized then $g_2$, that is adjusted with a naive idea of localization as a characteristic of a measure for a function support.    

\begin{table}[ht]
\centering
\caption{$UP_G$ and $UP_\lambda$: Example 1.}
\begin{tabular}{ccccccccc} 
\hline
 $f$ 
  & $\tilde{x}_0(f)$
  & $\tilde{t}_0(f)$
  & $V_{\lambda}(f)$
  & $V_{\lambda}(Ff)$
  & $UP_{\lambda}(f)$
  & $V_{G}(f)$
  & $V_{G}(Ff)$
  & $UP_{G}(f)$
	\\
   \hline
$f_1$ &  $[0,\,1/4)$ & $[0,\,4)$   &  $1/48$ & $16/3$    & $1/9$ & $1/28$ & $64/7$  & $16/49$  \\
$g_1$ &  $[3/4,\,1)$ & $[0,\,4)$   &  $1/48$ & $16/3$  & $1/9$ &  $1/28$ & $64/7$  &$16/49$  \\
$f_2$ &  $[0,\,1/8)$ & $[0,\,2)$   & $3/64$ & $8$  & $3/8$ & $4/21$ & $96/7$ & $128/49$ \\
$g_2$ &   $[3/4,\,7/8)$ & $[0,\,4)$   &  $71/64$ & $32/3$  & $71/6$ &$19/14$  & $255/14$& $4845/196$ \\
 \hline
\end{tabular}
\label{table1}
\end{table} 

{\bf Example 2.}
Here we discuss a dependence of a localization for a fixed function on a parameter $p$ of the Vilenkin group $G_p.$ Let us consider a function $f_1(x)= \mathbbm{1}_{\lambda^{-1}[0,\,1/4)}(x)$ and $p=2^k$,  $k\in \mathbb{N}$. We calculate $UP_G(f_1)$.

  (1)  If $k=1$, then $UP_G(f_1)=\displaystyle\frac{16}{49}$ (see Example 1.);

(2)  If $k=2$, then
$$
 V_{G}(f_1)=\frac{1}{\|f_1\|^2_{L_2(G)}}\min_{\tilde{x}} \int_{G}\|x\ominus\tilde{x}\|^2_{G}|f_1(x)|^2dx=
4 \min_{\widetilde{x}} \int_{\lambda^{-1}[0,\frac{1}{4})}\|x\ominus\tilde{x}\|^2_{G}dx
$$
$$
=4\int_{\lambda^{-1}[0,\frac{1}{4})}\|\tau\|^2_{G}d\tau=
4\sum_{i=1}^{\infty}\int_{\lambda^{-1}[\frac{1}{4^{i+1}},\frac{1}{4^i})}\|\tau\|^2_{G}d\tau=4\sum_{i=1}^{\infty}\left(\frac{1}{4^i}-\frac{1}{4^{i+1}}\right)4^{-2i}=\frac{1}{21}.
$$

$$
 V_{G}(Ff_1)=\frac{1}{\|Ff_1\|^2_{L_2(G)}}\min_{\tilde{x}} \int_{G}\|x\ominus\tilde{x}\|^2_{G}|Ff_1(x)|^2dx=\frac{1}{4}
\min_{\tilde{x}} \int_{\lambda^{-1}[0,4)}\|x\ominus\widetilde{x}\|^2_{G}dx
$$
$$
=\frac{1}{4}\int_{\lambda^{-1}[0,4)}\|\tau\|^2_{G}d\tau =
\frac{1}{4}\sum_{i=-1}^{\infty}\int_{\lambda^{-1}[\frac{1}{4^{i+1}},\frac{1}{4^i})}\|\tau\|^2_{G}d\tau=\frac{1}{4}
\sum_{i=-1}^{\infty}\left(\frac{1}{4^i}-\frac{1}{4^{i+1}}\right)4^{-2i}=\frac{256}{21}.
$$

Hence, $UP_G(f_1)=\displaystyle\frac{256}{441}$.

(3) If $k>2$, then

\[
 V_{G}(f_1)=\frac{1}{\|f_1\|^2_{L_2(G)}}\min_{\tilde{x}} \int_{G}\|x\ominus\tilde{x}\|^2_{G}|f_1(x)|^2dx=4 \min_{\widetilde{x}} \int_{\lambda^{-1}[0,\frac{1}{4})}\|x\ominus\tilde{x}\|^2_{G}dx\]\[
=4\int_{\lambda^{-1}\left[0,\frac{1}{2^k}\right)\oplus\left[\frac{1}{2^k},\frac{1}{4}\right)}\|\tau\|^2_{G}d\tau
=4\left(\sum_{i=1}^{\infty}\left(\frac{1}{(2^k)^i}-\frac{1}{(2^k)^{i+1}}\right)(2^k)^{-2i}+\left(\frac{1}{4}-\frac{1}{2^k}\right)\right)\]\[=
1-\frac{4}{2^k}+\frac{4}{2^k(2^{2k}+2^k+1)}.
\]

\[
 V_{G}(Ff_1)=\frac{1}{\|Ff_1\|^2_{L_2(G)}}\min_{\tilde{x}} \int_{G}\|x\ominus\tilde{x}\|^2_{G}|f_1(x)|^2dx=\frac{1}{4} \min_{\widetilde{x}} \int_{\lambda^{-1}[0,4)}\|x\ominus\tilde{x}\|^2_{G}dx\]\[
=\frac{1}{4}\int_{\lambda^{-1}\left[0,1\right)\oplus\left[1,4\right)}\|\tau\|^2_{G}d\tau
=\frac{1}{4}\left(\sum_{i=0}^{\infty}\left(\frac{1}{(2^k)^i}-\frac{1}{(2^k)^{i+1}}\right)(2^k)^{-2i}+(4-1)\cdot2^{2k}\right)\]\[=
\frac{3}{4}\cdot2^{2k}+\frac{1}{4}\cdot\frac{2^{2k}}{2^{2k}+2^k+1}.
\]

Therefore, $UP_G(f_1)=\displaystyle \left(1-\frac{4}{2^k}+\frac{4}{2^k(2^{2k}+2^k+1)}\right) \left(\frac{3}{4}\cdot2^{2k}+\frac{1}{4}\cdot\frac{2^{2k}}{2^{2k}+2^k+1}\right).
$

It is easy to see that time variance $V_G(f_1)$ goes to $1$, and frequency variance $V_G(Ff_1)$ goes to infinity as $k\to\infty.$


\section{Uncertainty product $UP_{G}$.}

In this section we concentrate on the uncertainty product corresponding to the metric $d_2.$ 
It turns out that the modified Gibbs derivative $\cal{D}$ plays a role of a usual derivative in this case. And since the Haar functions are the eigenfunctions of $\cal{D}$, it is possible to get representation for $UP_{G}$ using the Haar coefficients.  

\begin{teo} \label{UP_Haar}
Suppose $f\in L_2(G)\cap L_1(G),$ $\|\cdot\|_G f \in L_2(G)$, where ``dot'' $\cdot$   means the argument $x\in G$ of a function $f$,
and $\displaystyle f(x)=\sum_{\nu=1}^{p-1}\sum_{j\in\z}\sum_{k\in\z_+}c_{j,k}^{\nu}\psi_{j,k}^{\nu}(x)$. Then 
\be
\label{vF_haar}
\int_G \|t\|^2_G |Ff(t)|^2\,dt = \int_G  |{\cal D}f(t)|^2\,dt = \sum_{\nu=1}^{p-1}\sum_{j\in\z} \sum_{k\in\z_+} |p^j c_{j,k}^{\nu}|^2
\ee  
\be
\label{v_haar}
\int_G \|x\|^2_G |f(x)|^2\,dx = \int_G  |{\cal D}Ff(x)|^2\,dx =  \sum_{\nu=1}^{p-1} \sum_{j\in\z} \sum_{k\in\z_+} |p^j d_{j,k}^{\nu}|^2,
\ee  
where $d_{j,k}^{\nu}$, $j\in\z,$ $k\in\z_+$, $\nu=1,\dots,p-1$, are the coefficients in the Haar series for the function $Ff,$ that is $\displaystyle Ff(t)=\sum_{\nu=1}^{p-1}\sum_{j\in\z}\sum_{k\in\z_+}d_{j,k}^{\nu}\psi_{j,k}^{\nu}(t).$
\end{teo}

\textbf{Proof.} 
By the definition of the modified Gibbs derivative and the Plancherel equality we get
$$
\int_G \|t\|^2_G |Ff(t)|^2\,dt = \int_G  |F{\cal D}f(t)|^2\,dt= 
\int_G  |{\cal D}f(t)|^2\,dt.
$$
Expanding a function in the Haar series and applying Corollary \ref{eigen_Haar}, we get
$$
\int_G  |{\cal D}f(t)|^2\,dt 
=
\int_G  \left|\sum_{\nu=1}^{p-1}\sum_{j\in\z}\sum_{k\in\z_+}c^{\nu}_{j,k}{\cal D}\psi^{\nu}_{j,k}(t)\right|^2\,dt
$$
$$
=
\int_G  \left|\sum_{\nu=1}^{p-1}\sum_{j\in\z}\sum_{k\in\z_+}c^{\nu}_{j,k}p^j\psi^{\nu}_{j,k}(t)\right|^2\,dt
= \sum_{\nu=1}^{p-1}\sum_{j\in\z} \sum_{k\in\z_+} |p^j c_{j,k}^{\nu}|^2
$$
The last equality follows from the orthonormality of the Haar system. Equality (\ref{v_haar}) is proved analogously to (\ref{vF_haar}).
 \hfill $\Box$

\begin{rem}
Formally, it is possible to write 
$
\int_G \lambda^2(x) |Ff(x)|^2 \, dx = \int_G |f^{[1]}(x)|^2\, dx 
$
and to try to represent $UC_{\lambda}$ in terms of eigenfunctions of the Gibbs derivative $f^{[1]}$ in the case of the Cantor group. (The Gibbs derivative is defined for functions defined on the Cantor group only.) However, the Gibbs differentiation is not a local operation, that is 
$(f \mathbbm{1}_E)^{[1]}\neq f^{[1]} \mathbbm{1}_E$, see also discussion in \cite{LS15}. So, usage of Walsh functions instead of Haar basis might give interesting results for periodic functions only.    
\end{rem}

We did not found in the  literature a formula expressing $d_{j,k}^{\mu}$ in terms of $c_{j,k}^{\nu}$. So we obtain this formula in the following lemma.  
 \begin{lem}
\label{lem_d_c} 
  Suppose $f\in L_2(G)$ and the coefficients $c_{j,k}^{\nu}$, $d_{j,k}^{\mu}$, $j\in\z,$ $k\in\z_+$, $\nu, \mu=1,\dots,p-1$,  are defined in Theorem \ref{UP_Haar}. Then 
\begin{gather}
d_{j,k}^{\mu} = 
\sum_{\nu=1}^{p-1} p^{q_0/2} b^{\nu}_{k} 
+ p^{j/2} \sum_{\nu=1}^{p-1} c^{\nu}_{-j-1,0} {\rm exp}\left(-\frac{2\pi i \nu \mu}{p}\right) \delta_{k,0} + p^{j/2} \sum_{i=-\infty}^{-j-2} \sum_{\nu=1}^{p-1} c^{\nu}_{i,0} \delta_{k,0},
\label{d_c}
\end{gather}
 where $ \displaystyle b^{\nu}_k = p^{-q_0} \sum_{n=0}^{p^{q_0}-1} c^{\nu}_{q_0-j,n+(p-\mu)p^{q_0}} 
\chi(\lambda^{-1}(n),D^{-q_0}\lambda^{-1}(k))$ is the $k$-th term of the  discrete Vilenkin-Chrestenson
transform of $(c_{q_0-j,n+(p-\mu)p^q_0}^{\nu})_{n=0}^{p^q_0-1},$ $\displaystyle q_0=\left[\log_p \frac{k}{p-\nu}\right]$, 
and $\delta_{0,0}=1$, and $\delta_{k,0}=0,$ if $k\neq 0.$ 
\end{lem}

\textbf{Proof.} 
Using the Plancherel equality and (\ref{hat_haar}), we get
$$
d^{\mu}_{j,k}=\int_{G} Ff(x)\overline{\psi^{\mu}_{j,k}(x)}\,dx = 
\int_{G} f(x)\overline{F\psi^{\mu}_{j,k}(x)}\,dx =
\sum_{\nu=1}^{p-1} \sum_{i\in\z} \sum_{n\in\z_+}c^{\nu}_{i,n}\int_{G} \psi^{\nu}_{i,n}(x)\overline{F\psi^{\mu}_{j,k}(x)}\,dx
$$
$$
 = \sum_{\nu=1}^{p-1} \sum_{i\in\z} \sum_{n\in\z_+}c^{\nu}_{i,n}\int_{G} \psi^{\nu}_{i,n}(x) 
p^{-j/2} \overline{\chi(\lambda^{-1}(k),D^{-j} x)} \mathbbm{1}_{I_{-j}\oplus \lambda^{-1}((p-\mu)p^j)}(x) \,dx.
$$
Since $\supp \psi^{\nu}_{i,n}= \lambda^{-1}([n p^{-i},(n+1)p^{-i})),$ it follows that the last expression takes the form
$$
\sum_{\nu=1}^{p-1} \sum_{i=-j}^{\infty} \sum_{n=(p-\mu)p^{i+j}}^{(p-\mu+1)p^{i+j}-1} c^{\nu}_{i,n}\int_{G} \psi^{\nu}_{i,n}(x) 
p^{-j/2} \overline{\chi(\lambda^{-1}(k),D^{-j} x)}  \,dx
$$
$$
+ 
p^{-j/2} \left(\sum_{\nu=1}^{p-1} c^{\nu}_{-j-1,0} {\rm exp}\left(-\frac{2\pi i\nu \mu}{p}\right)+
\sum_{\nu=1}^{p-1} \sum_{i=-\infty}^{-j-2} c^{\nu}_{i,0} 
\right)
$$
$$
\times
\int_{G} \overline{\chi(\lambda^{-1}(k),D^{-j} x)} \mathbbm{1}_{I_{-j}\oplus \lambda^{-1}((p-\mu)p^j)}(x) \,dx =:
S_1+S_2.
$$
For the first sum by (\ref{hat_haar}) we note that 
$$
\int_{G} \psi^{\nu}_{i,n}(x) 
 \overline{\chi(\lambda^{-1}(k),D^{-j} x)}  \,dx= F \psi^{\nu}_{i,n}(D^{-j}\lambda^{-1}(k))
$$
$$
= p^{-i/2} \chi(n,D^{-i-j}\lambda^{-1}(k)) 
\mathbbm{1}_{I_{-i}\oplus\lambda^{-1}((p-\nu)p^i)}(D^{-j}\lambda^{-1}(k)).
$$
Therefore, the first sum takes the form
$$
S_1 = \sum_{\nu=1}^{p-1} \sum_{i=-j}^{\infty} \sum_{n=(p-\mu)p^{i+j}}^{(p-\mu+1)p^{i+j}-1}p^{-(j+i)/2} c^{\nu}_{i,n} \chi(\lambda^{-1}(n),D^{-i-j}\lambda^{-1}(k))
\mathbbm{1}_{I_{-i-j}\oplus\lambda^{-1}((p-\nu)p^(i+j))}(\lambda^{-1}(k))
$$
$$
= \sum_{\nu=1}^{p-1} \sum_{q=0}^{\infty} p^{-q/2} \sum_{n=0}^{p^{q}-1} c^{\nu}_{q-j,n+(p-\mu)p^{q}} 
\chi(\lambda^{-1}(n),D^{-q}\lambda^{-1}(k))
\mathbbm{1}_{I_{-q}\oplus\lambda^{-1}((p-\nu)p^q)}(\lambda^{-1}(k)).
$$
Since $\mathbbm{1}_{I_{-q}\oplus\lambda^{-1}((p-\nu)p^q)}(\lambda^{-1}(k))=1$ for 
$(p-\nu)p^q\le k < (p-\nu+1)p^q$ and  $\mathbbm{1}_{I_{-q}\oplus\lambda^{-1}((p-\nu)p^q)}(\lambda^{-1}(k))=0$ for the remaining $k$, 
and since the inequality $(p-\nu)p^q\le k < (p-\nu+1)p^q,$ $q\in \mathbb{Z}_+$ is equivalent to 
$\displaystyle q=\left[\log_p \frac{k}{p-\nu}\right]$, it follows that the only nonzero term in the sum 
$\sum_{q=0}^{\infty}$ has the number $\displaystyle q_0:=\left[\log_p \frac{k}{p-\nu}\right]$. 
So
$$
S_1 = \sum_{\nu=1}^{p-1} p^{-q_0/2} \sum_{n=0}^{p^{q_0}-1} c^{\nu}_{q_0-j,n+(p-\mu)p^{q_0}} 
\chi(\lambda^{-1}(n),D^{-q_0}\lambda^{-1}(k)).
$$
By (\ref{71.1}) we notice that up to the multiplication by a constant  the inner sum in the last expression is the $k$-th term of   
the  discrete Vilenkin-Chrestenson
transform of the vector $(c_{q_0-j,n+(p-\mu)p^q_0}^{\nu})_{n=0}^{p^q_0-1}.$ Denote this term by $b_k^{\nu}$. Finally, for $S_1$ we get 
$$
S_1(x)= \sum_{\nu=1}^{p-1} p^{q_0/2} b_k^{\nu}.
$$
  Thus, the first sum takes the desired form.
To conclude the proof it remains to calculate the following part of the second sum  
$$
\int_{G} \overline{\chi(k,D^{-j} x)} \mathbbm{1}_{I_{-j}\oplus \lambda^{-1}((p-\mu)p^j)}(x) \,dx =
p^j \int_{G} \overline{\chi(k, x)} \mathbbm{1}_{I\oplus \lambda^{-1}(p-\mu)}(x) \,dx
$$
$$ =
 p^j \int_{I} \overline{\chi(k, x\ominus \lambda^{-1}(p-\mu))}\,dx = 
p^j \int_{I} \overline{\chi(k, x)}\,dx = p^j \delta_{k,0},
$$
where $\delta_{0,0}=1$, and $\delta_{k,0}=0,$ if $k\neq 0.$
\hfill $\Box$

It is easy to see from (\ref{vF_haar}) that $\min \int_G \|t\|^2_G |Ff(t)|^2\,dt=0$ and $\max \int_G \|t\|^2_G |Ff(t)|^2\,dt=\infty$ under the restriction $\|f\|_{L_2(G)}=1.$

Formulas (\ref{vF_haar}) and (\ref{v_haar}) allow for the following result on estimation of Fourier-Haar coefficients for functions defined on the Vilenkin group.

\begin{cor} \label{Haar_coef}
Suppose  $\|\cdot\|_G Ff \in L_2(G)$,
 and $\displaystyle f(x)=\sum_{\nu=1}^{p-1}\sum_{j\in\z}\sum_{k\in\z_+}c_{j,k}^{\nu}\psi_{j,k}^{\nu}(x)$. Then the series 
$
\displaystyle \sum_{\nu=1}^{p-1}\sum_{j\in\z} \sum_{k\in\z_+} |p^j c_{j,k}^{\nu}|^2
$  
is convergent.
\end{cor}

\section*{Acknowledgments}
The authors are supported by Volkswagen Foundation.
The second author is supported by the RFBR, grant \#15-01-05796, 
and  by Saint Petersburg State University, grant  \#9.38.198.2015.




\begin{thebibliography}{0}



\bibitem{B} E.~Breitenberger,
 Uncertainty measures and uncertainty relations for angle observables,
Found. Phys.    \textbf{15} (1985), 353--364.

	
\bibitem{G} B.I. Golubov, Elements of dyadic analysis, [in Russian], Moscow, LKI, 2007.


\bibitem{GES} B.I. Golubov, A.V. Efimov, and V.A. Skvortsov, Walsh series and transforms, English transl.: Kluwer, Dordrecht, 1991.

\bibitem{Erb} W. Erb, Uncertainty principles on compact Riemannian manifolds, Appl. Comput. Harmon. Anal., \textbf{29} (2010), 182-197.

\bibitem{F} Yu. A. Farkov, Multiresolution analysis and wavelets
on Vilenkin groups, FACTA UNIVERSITATIS (NIS) 
SER.: ELEC. ENERG. vol. 21, no. 3, December 2008, 309-325.

\bibitem{H} W.~Heisenberg, The actual concept of quantum theoretical kinematics and mechanics,    
Physikalische Z.  \textbf{43}  (1927), 172.

\bibitem{HR}     E. Hewitt, K. A. Ross, Abstract Harmonic Analysis.
Springer-Verlag, New York, 1963, 1979.

\bibitem{27}  A.~V.~Krivoshein,  E.~A.~Lebedeva,  Uncertainty Principle for the Cantor Dyadic Group,  J. Math. Anal. Appl, \textbf{42} (2015),  1231-1242.

\bibitem{Lang}   W. C. Lang. Orthogonal wavelets on the Cantor dyadic group, SIAM J. Math. Anal.~--- 1996.~--- Vol.~27.~--- P.~305--312.

\bibitem{LS15} E.~Lebedeva, M.~Skopina. Walsh and wavelet methods for differential equations on the Cantor group // J. Math. Anal.  Appl.~--- 2015.~--- Vol.~430.~--- No.~2.~---  P.~593--613.


 \bibitem{PrSi} J. F.Price, A. Sitaram, Local uncertainty inequalities for locally compact groups,
 Trans. of AMS, \textbf{308} 1 (1988), 105--114. 

	\bibitem{Sch} E. Schr\"odinger, About Heisenberg uncertainty relation, Proc. of The Prussian Acad. of Scien. XIX (1930)
296--303.


\bibitem{SWS}F.~Schipp, W.~R.~Wade, P.~Simon. Walsh series. 
An introduction to dyadic harmonic analysis.~---  Academiai Kiado, Budapest, 1990.


\end{thebibliography}
\end{document}